\documentclass[11pt, reqno]{amsart}
\newtheorem{thm}{Theorem}
\newtheorem{lem}[thm]{Lemma}
\newtheorem{prop}[thm]{Proposition}
\newtheorem{cor}[thm]{Corollary}
\theoremstyle{remark}

\theoremstyle{definition}
\newtheorem{defi}{Definition}
\newcommand{\R}{\mathbb{R}}

\newcommand{\N}{\mathbb{N}}
\newcommand{\Q}{\mathbb{Q}}

\newcommand{\vect}{\mathit{Vect}}

\newcommand{\del}{\partial}

\newcommand{\alg}[2]{\mathit{{#1}}({#2})}

\newcommand{\id}{\mathrm{id}}
\newcommand{\algproj}{\mathit{sl}_{n+1}}
\newcommand{\mdiff}[3]{\mathcal{M}_{\mathit{diff}}^{#1}({#2},{#3})}
\newcommand{\msymb}[2]{\mathcal{S}^{(#1)}_{#2}}

\newcommand{\casi}{\mathcal{C}^\mathit{op}}
\newcommand{\casisymb}{\mathcal{C}^t}
\newcommand{\casinil}{N_\mathcal{C}}
\newcommand{\Lie}[1]{L_{{#1}}}
\renewcommand{\S}[1]{\mathcal{S}_{{#1}}}
\newcommand{\hwv}[1][k,l,q]{\vec{v}_\mathit{H}({#1})}
\newcommand{\pent}[1]{\lfloor {#1} \rfloor}

\title[Projective symbols for bidifferential operators]{Projectively equivariant symbol calculus for bidifferential 
operators}
\author{Fabien Boniver}
\address{
 Universit\'e de Li\`ege, Institut de Math\'ematique, B37\\ Grande Traverse, 
 12, B-4000 Sart Tilman (Li\`ege), Belgium
}
\email{f.boniver@ulg.ac.be}
\begin{document}
\begin{abstract}
    We prove the existence and uniqueness of a \emph{projectively 
    equivariant symbol map}, which is an isomorphism between 
    the space of bidifferential operators acting on tensor densities over 
    $\R^n$ and that of their symbols, when both are 
    considered as  modules over an imbedding of $\alg{sl}{n+1,\R}$ into 
    polynomial vector fields.
    
    The coefficients of the bidifferential operators are densities of 
    an arbitrary weight.  We obtain the result for all values of this 
    weight, except for a set of \emph{critical} ones, which does not 
    contain $0$.  In the case of second order operators, we give explicit formulas and examine in detail the 
    critical values.
\end{abstract}
\maketitle

\section{Introduction}
The spaces of linear differential operators acting on tensor densities 
over a smooth manifold $M$ naturally constitute representations of 
its Lie algebra of vector fields.  Let $\mathcal{D}_{\lambda}(M)$ 
denote the module of differential operators acting on densities of 
weight $\lambda$.  The classification of these modules --- regarding $\lambda$ as a 
parameter --- was recently obtained (see~\cite{do,go,lmt}).

A new and fruitful approach to this classification consists in comparing the 
action of vector fields on differential operators and on their 
symbols.  This method is exposed and applied to the classification of 
quotient modules in~\cite{lo}.  One proves the existence and takes 
advantage of the properties of a \emph{projectively equivariant 
symbol map}
\[
\sigma_{\lambda}:\mathcal{D}_{\lambda}(\R^n)\rightarrow 
\bigoplus_{k\geq 0}\Gamma(\vee^k T \R^n)
\]
which has the following properties.  It is the unique isomorphism of 
representations of the \emph{projective Lie algebra} to preserve the 
principal symbol of its arguments.  The latter algebra is the 
subalgebra of polynomial vector fields generated by
\begin{equation}\label{eq.sl}
\frac{d}{dx^i}, x^j \frac{d}{dx^i},\text{ and } 
x^j\sum_{i}\frac{d}{dx^i},\quad i,j\in\{1,\ldots,n\}.
\end{equation}
We shall denote it by $\algproj$ since it is isomorphic to 
$\alg{sl}{n+1,\R}$ (see for instance~\cite[p. 6]{fuks}).

It is worth noticing that the inverse function of $\sigma_{\lambda}$ 
maps functions on the cotangent bundle of the Euclidean space that 
are polynomial on the fibre --- which are nothing but symmetric 
contravariant tensor fields --- onto differential operators.  It may 
therefore be viewed as a quantization procedure (cf.~\cite{do} and 
references therein)
that preserves the infinitesimal symmetries of 
the base space.  It may also be used to define equivariant star-products.
Two further steps are made in~\cite{dlo}. On the one 
hand, the subalgebra prescribing the infinitesimal invariance is there 
changed to be made of conformal transformations; on the other hand, a 
wider class of differential operators with densities of an arbitrary 
weight as coefficients is studied.

This paper is intended to be a first step towards the classification of 
modules of bilinear \emph{bidifferential} operators.  It is clearly 
a natural approach to wonder about the existence of an equivariant symbol 
map for these.
However, computing an analogue of $\sigma_{\lambda}$ using the 
prescription of equivariance turns out to be more intricate than 
in the linear case.  Therefore, we have used techniques from the 
``conformal'' case~(cf.~\cite{dlo}) and collected interesting 
intermediate results about the structure of the considered modules, 
such as their Casimir operators and the spectrum of these.

These aspects lead us towards a more algebraic study of modules of 
bidifferential operators over tensor densities.
It is  worth noticing that these have already been studied, see for 
instance~\cite{gro}.

Our paper is organized as follows.  In Section~\ref{sec.defi}, we recall 
basic definitions and 
choose notations. In Section~\ref{sec.casi}, we first compute 
the Casimir operators of $\algproj$ acting on multi\-differential 
operators and tensor fields, and then determine the spectrum of the 
latter when specialized to bidifferential operators.
We show in Section~\ref{sec.bidiffl} that the former has the 
same spectrum, which allows to define an equivariant symbol map.  
This is possible due to a strong assumption on the weight of the 
coefficients of the operators, which we shall call \emph{shift}.  
We next weaken this assumption and give a lower bound function 
for the remaining \emph{critical} shift values.  Section~\ref{sec.secord} is 
devoted to an explicit computation of the symbol map restricted to 
second order operators.  This example demonstrates that when the shift 
is critical, a
symbol map may either not exist or not be unique.
\section{Definitions}
\label{sec.defi}
\subsection{Densities}
Let $\lambda$ be a real number.
A \emph{density of weight $\lambda$} on $\R^n$ is a function 
$\phi:\Lambda^n T\R^n\setminus \{0\}\rightarrow\R$ such that
\[
\phi(t\omega)=\vert t \vert^\lambda \phi(\omega),\quad \forall 
t\in\R\setminus 0,\forall \omega\in\Lambda^n T\R^n\setminus \{0\}.
\]
Let $E^\lambda$ denote the one-dimensional vector space of 
densities of weight $\lambda$ on $\R^n$.  

For the sake of simplicity, we shall also call densities the sections of 
$\Delta^\lambda(\R^n)=\R^n\times E^\lambda$ viewed as a fibre bundle over $\R^n$.

The pull-back of $\phi$ along the 
flow of a vector field $X$ allows to define the Lie derivative of 
$\phi$ in the direction of $X$, yielding the expression~:
\begin{equation}\label{eq:Liedens}
\Lie{X}\phi=\sum_{i}X^i\del_{i}\phi+\lambda \sum_{i} \del_{i} X^i \phi.
\end{equation}
We shall write $\del_{i}$ for both the derivative and the unit vector along the $i$th axis.

\subsection{Multidifferential operators}
Let
$
\mdiff p \lambda \mu
$
be the space of multilinear multi\-differential operators transforming $p$ 
densities of weights \linebreak
$\lambda_1,\ldots,\lambda_p$ ($\lambda$ is the multi-index of these 
weights) into one
density of weight $\mu$.

Such an operator is said to have \emph{(total) order} $r$ if it  may be written
\begin{equation}\label{eq.op}
(f_{1},\ldots,f_{p})\mapsto\sum 
A_{\alpha^1,\ldots,\alpha^p}\,D^{\alpha^1}\!\! f_{1}\cdots\, D^{\alpha^p} 
\!\!f_{p}
\end{equation}
with ${\sum_{j=1}^p \vert\alpha^j\vert\leq r}$ , where the $\alpha^{j}$ 
are multi-indices and 
\[
D^{\alpha^j}=\del_{\alpha^j_{1}}\cdots\del_{\alpha^j_{k}},\quad 
\forall j, \quad(\vert\alpha^j\vert=k).
\]
The coefficients $A_{\alpha^1,\ldots,\alpha^p}$ in~(\ref{eq.op}) are densities 
of weight 
$
\delta =\mu-\sum_{i}\lambda_{i}.
$
\emph{Throughout this paper, unless otherwise stated, we shall 
consider $\delta,\lambda$ and $\mu$ as above.  We shall moreover refer 
to $\delta$ as the \emph{shift} value.}

The space of those operators is a module over the Lie 
algebra $\vect(\R^n)$ of vector fields 
over $\R^n$. It is filtered by the order.  Indeed, the 
Lie derivative of an $r$th order operator $T\in\mdiff p \lambda \mu$ is given by
\begin{equation}\label{eq.liederop}
\Lie{X}T=\Lie{X}\circ T-\sum_{i=1}^p 
T(\ldots,\overbrace{\Lie{X}\cdot}^{\mbox{$i$th arg.}},\ldots),
\end{equation}
which is again an $r$th order operator.  
The Lie derivatives of the right hand side of the latter formula are those 
defined by~$(\ref{eq:Liedens})$ with the suitable weight of densities.

\subsection{Symbols}
The polynomial functions to be associated to those operators, which 
we also name \emph{symbols}, are the elements of
\begin{equation}\label{eq:symbolspace}
\msymb p {\delta}=\bigoplus_{i_{1},\ldots,i_{p}\geq 
0}\Gamma(\vee^{i_{1}}(T\R^n)\otimes\cdots\otimes\vee^{i_{p}}(T\R^n)\otimes\Delta^{\delta}(\R^n)),
\end{equation}
where the final term takes the shift value into account.
The Leibniz rule and the Lie derivatives of vector fields and 
densities endow this space with a 
natural structure of $\vect(\R^n)$-module.

\subsection{Projectively equivariant symbol map}
The association of the tensor field
\begin{equation}\label{eq.symb}
(\xi^1,\ldots,\xi^p)\mapsto A_{\alpha^1,\ldots,\alpha^p} 
(\xi^1)^{\alpha^{1}}\cdots (\xi^p)^{\alpha^p} \in\msymb{p}{\delta}
\end{equation}
to the operator defined in~(\ref{eq.op}) is an isomorphism of vector 
spaces between $\msymb{p}{\delta}$ and $\mdiff p \lambda \mu$.  It is 
well known, and easy to check, that it is not equivariant with 
respect to $\algproj$.  Remember that we defined  in~(\ref{eq.sl}) this subalgebra of 
vector fields, with respect to which we want to impose equivariance.

The \emph{projectively equivariant symbol map} we are looking for is a linear 
bijection 
\[
\sigma_{\lambda,\mu}:\mdiff 2 \lambda \mu\rightarrow \msymb{2}{\delta}
\]
that preserves the principal symbol of its arguments and such that 
\[
\sigma_{\lambda,\mu}(\Lie{X} T)=\Lie{X}(\sigma_{\lambda,\mu} T)
\]
for all $T\in\mdiff 2 \lambda\mu$ and $X\in\algproj$.

\subsection{Polynomial form}
Though it is not equivariant, the symbolic notation~(\ref{eq.symb}) of 
a multidifferential operator is very useful to handle computations 
with operators that might otherwise become tedious.  Moreover, one can 
take advantage of the following rule, ``\`a la Fourier''.  For the 
sake of simplicity, assume that $D$ is a differential operator, whose symbolic 
notation is the polynomial function $\xi\mapsto P(\xi)$, in which $\xi$
represents the derivatives in $D$ that affect 
its argument, say $f$.  Then, if $X$ is a vector field, the operator 
\[
f\mapsto D(\sum X^i \del_{i} f)
\]
admits the polynomial function $\xi\mapsto P(\xi+\eta) <X,\xi>$ as 
its symbolic notation, if one agrees on denoting by $\eta$ the 
derivatives affecting $X$.  This provides not only a convenient notation for the 
action of $D$ and $X$, but also an easy way to keep track of the applications 
of derivatives on products.  Another example is the symbolic form of 
the Lie derivative of a multidifferential operator~(\ref{eq.LieDerT}) 
below.

\section{Casimir operators}
\label{sec.casi}
In this section, we shall compute the Casimir operators $\casi$ of $\mdiff p 
\lambda \mu$ and $\casisymb$ of $\msymb p \delta$ as representations of the 
projective Lie algebra.  

Let the derivatives of $T\in\mdiff p \lambda \mu$ be represented by 
the covariant
variables $\xi^1,\ldots,\xi^p$.  Define
\[
\tau_{i,\zeta} T = T(\ldots,\overbrace{\cdot+\zeta}^{\mbox{i-th arg.}},\ldots)-T.
\]
If $\eta$ denotes the derivatives affecting the coefficients of $T$ 
and if $\zeta$ is similarly related to a vector field $X$, 
then one may write the symbolic form of the Lie derivative of $T$ in 
the direction of $X$:
\begin{multline}\label{eq.LieDerT}
\Lie{X} T =\\
<X,\eta> T - \sum_{i=1}^p (<X,\xi^i>+\lambda_i 
<X,\zeta>)\tau_{i,\zeta} T
+ \delta <X,\zeta> T.
\end{multline}
This results of a direct computation from the 
definition~(\ref{eq.liederop}) of the Lie derivative 
on $\mdiff p\lambda \mu$.

If $P$ is the polynomial form of $T$ then
\begin{equation}\label{eq.LieDerP}
\Lie{X} P =
<X,\eta> P - \sum_{i=1}^p <X,\xi^i>(\zeta D_{\xi^i}) P
+ \delta <X,\zeta> P.
\end{equation}

To achieve the computation of the Casimir operators, we fix two bases of $\algproj$ dual to 
each other with respect to some non-degenerate symmetric bilinear 
form (and therefore a multiple of the Killing form of the 
algebra).
  We shall consider the 
 vector fields~:
\begin{equation}\label{eq.basis}
\begin{array}{llllllll}
    e_{ij} & = &-x^j\del_{i},& \quad (i\neq j) & e_{ij}^*& = & e_{ji}, 
    &\quad (i\neq 
    j) \\
    e_{ii} &=& -x^i\del_{i}-\sum_{k}x^k \del_{k}& & e_{ii}^* &=& -x^i\del_{i}& \\
    e_{i} &=& -\del_{i} && e_{i}^* &=& \epsilon^i &\\
    \epsilon^i &=& x^i\sum_{k} x^k\del_{k}& & \epsilon^{i*} &=& e_{i}&
\end{array}
\end{equation}
for $i,j\in\{1,\ldots,n\}$, where $e_{ij}$ and $e_{ij}^*$ are dual, 
and so on.

\begin{prop}
The Casimir operator 
\begin{equation}\label{eq.CasiDef}
\casi=\sum_{i\neq j}L_{e_{ij}}\circ L_{e_{ij}^*}+\sum_{i}(L_{e_{ii}}\circ 
L_{e_{ii}^*}+L_{e_{i}}\circ L_{e_{i}^*}+L_{\epsilon^i}\circ 
L_{\epsilon^{i*}})
\end{equation}
of $\algproj$ acting on $\mdiff p \lambda \mu$ 
equals the sum
\[
\casi=\casisymb+\casinil
\]
where
\begin{multline*}
\casisymb=n(n+1)\delta(\delta - 1) \id + 2((n+1)(1-\delta))\sum_k E_{\xi^k} \\
\mbox{}
+\sum_{k,l} \sum_{i,j} \xi^k_i \xi^l_j D_{\xi^k_j} 
D_{\xi^l_i}  
 + \sum_{k,l} \sum_{i,j} \xi^k_i \xi^l_j D_{\xi^k_i} 
D_{\xi^l_j}
\end{multline*}
is the Casimir operator of $\algproj$ acting on $\msymb p \delta$
and
\[
\casinil=2 \sum_k (E_{\xi^k} + (n+1) \lambda_k) (\eta D_{\xi^k}),
\]
if $\eta$ denotes the partial derivatives of the 
coefficients of an argument of $\casi$ and 
$E_{\xi^k}=\sum_{i}\xi^k_{i}D_{\xi^k_{i}}$.
\end{prop}
\begin{proof}
Let $T\in\mdiff p \lambda \mu$.
    We first compute the Lie derivatives of T  in the directions of the 
vector fields given by~(\ref{eq.basis}).  One gets
\[
\begin{array}{lcl}
L_{e_{ij}} T &=& -x^j\del_i T + \sum_{k=1}^p \xi^k_i D_{\xi^k_j} T 
\\[1ex]
L_{e_{ii}} T &=& -x^i\del_i T - E.T +\sum_k (E_{\xi^k} T + \xi^k_i 
D_{\xi^k_i} T -
\delta(n+1) T)\\[1ex]
L_{e_i} T &=& -\del_i T \\[1ex]
L_{\epsilon^i} T&=& x^i E.T+ \delta (n+1) x^i T \\[1ex] 
&&\mbox{}- \sum_k (\xi^k(x) D_{\xi^k_i} T + x^i E_{\xi^k} T +
E_{\xi^k} D_{\xi^k_i} T + (n+1) \lambda_k D_{\xi^k_i} T) \\[1ex]
L_{e^*_{ii}} T&=&-x^i \del_i T + \sum_k \xi^k_i D_{\xi ^k _i} T - 
\delta T
\end{array}
\]
where $E = \sum_i x^i\del_i$.

The Casimir operator intertwines $\mdiff p \lambda \mu$ with itself 
as a module over $\algproj$. 
It is in particular invariant with respect to the constant vector fields and has 
therefore constant coefficients.   To compute the Casimir operator, 
we thus only need to collect and sum  terms with constant coefficients  
in~(\ref{eq.CasiDef}).  Applying the formulas for Lie derivatives once 
again, we get
\begin{equation}\label{eq:Casiparts}
\begin{array}{llcl}
\mbox{for} & L_{e_{ij}}\circ L_{e_{ij}^*} T &:& \sum_{i\neq j} (\sum_k \xi^k_i 
D_{\xi^k_i} T +
\sum_{k,l} \xi^k_i \xi^l_j D_{\xi^k_j} D_{\xi^l_i} T), \\[1ex]
\mbox{for} & L_{e_{ii}}\circ L_{e_{ii}^*} T &:& (2-2(n+1)\delta) \sum_k E_{\xi^k} 
T \\[1ex]
&&&\mbox{}+ \sum_{i,j}
\sum_{k,l} \xi^k_i \xi^l_j D_{\xi^k_i} D_{\xi^l_j} T   \\[1ex]
&&&\mbox{}+\sum_i \sum_{k,l} \xi^k_i
\xi^l_i D_{\xi^k_i} D_{\xi^l_i} T + n(n+1) \delta^2 T,\\[1ex]
\mbox{for} & L_{e_i}\circ L_{\epsilon^i} T &:&  \sum_k 
(E_{\xi^k} +(n+1)
\lambda_k) (\eta D_{\xi^k}) T\\[1ex]
&&&\mbox{}- n(n+1)\delta T +(n+1) \sum_k E_{\xi^k} T \\[1ex]
\mbox{and for} & L_{\epsilon^i}\circ L_{e_i} T &:& \sum_k 
(E_{\xi^k}+(n+1)\lambda_k)(\eta D_{\xi^k}) T.
\end{array}
\end{equation}
 Summing these terms  gives the announced 
expression.
\end{proof}

\subsection{Spectrum of $\casisymb$}
We shall now deal with bidifferential operators.  Instead of $\xi^1$ and 
$\xi^2$, we shall write $\alpha$ 
and $\beta$ the
covariant variables symbolizing their partial derivatives.

From now on, we shall also assume that the dimension $n$ of 
the Euclidean space is at least $2$. One will find a discussion 
of the case $n=1$ in Section~\ref{sec.n1}.

We are now going to split $\msymb{2}{\delta}$ into spaces of homogeneous 
eigenvectors of $\casisymb$.  
Observe first that the Lie derivatives in direction of linear vector fields only contribute 
to the symbolic part of $\casi$.  Moreover, for any 
such vector field $X$, 
\[
\Lie{X}\circ\Lie{X^*} = \rho(DX)\circ \rho(DX^*),
\]
where $\rho$ denotes the natural representation of $\alg{gl}{n,\R}$ 
on $\msymb{2}{\delta}$.
Together with expression~(\ref{eq:Casiparts}), this equality shows that
on each space 
$\vee^k T\R^n\otimes\vee^l T\R^n\otimes\Delta^\delta(\R^n)$,  $\casisymb$ can be written
\[
\mathcal{C}_{\alg{gl}{n,\R}}+ r \,\id \quad (r\in\N),
\]
where $\mathcal{C}_{\alg{gl}{n,\R}}$ denotes the Casimir operator 
associated to the representation $\rho$ above on the space $\vee^k 
T\R^n\otimes\vee^l T\R^n$. 
Its irreducible $\alg{gl}{n,\R}$-submodules will 
therefore be made up of eigenvectors of $\casisymb$.

It is well known that, for  $k,l\in \N$,
\[
\vee^k T\R^n\otimes \vee^l T\R^n=\bigoplus_{q=0}^{\mathrm{min}(k,l)} 
V_{k,l,q}
\]
where $V_{k,l,q}$ is an irreducible $\alg{gl}{n,\R}$-submodule 
described by a Young ta\-bleau with two lines of respective lengths 
$k+l-q$ and $q$.  It is easy to check that  
\begin{equation}\label{eq:hwv}
\hwv : \alpha\otimes\beta\in\R^{n*}\otimes \R^{n*}\mapsto
(\alpha_1 \beta_2 - \alpha_2 \beta_1)^{q} \alpha_1^{k-q} \beta_1^{l-q}
\end{equation}
is a highest weight vector of $V_{k,l,q}$.
We shall denote the space
\[
\Gamma(V_{k,l,q}\otimes\Delta^\delta\R^n)
\]
by $\S{k,l,q}$, whenever $0\leq q\leq \min(k,l)$, and take this to be 
$\{0\}$ otherwise. 
 
\begin{lem}\label{lem.eigen}
The restriction of $\casisymb$ to $\S{k,l,q}$ equals
\[
n(n+1)\delta(\delta-1)-2((n+1)\delta -n+q)(k+l)+2(k+l)^2+2q(q-1)
\]
times the identity.
\end{lem}
\begin{proof}
Just evaluate $\casisymb$ on $\hwv$.
\end{proof}

Depending only on the shape of the Young tableau, this value contains $k$ and $l$ through their sum.  
Let us denote it by 
$\gamma_{k+l,q}$.  Notice that $\gamma_{i,p}=\gamma_{i,p'}$ only if 
$p=p'$.
We can thus gather eigenspaces of same eigenvalue and total tensor 
degree.  We define
\[
\S{(i,p)}=\bigoplus_{k=p}^{i-p} \S{k,i-k,p}.
\]

We have just shown how each element of $\msymb{2}{\delta}$ writes as a sum of homogeneous eigenvectors of 
$\casisymb$~:
\[
\msymb{2}{\delta}=\bigoplus_{i\geq 0}\,\bigoplus_{0\leq p\leq \pent{i/2}} 
\S{(i,p)},
\]
if $\pent{x}$ is the highest natural number less than or equal to $x$.
\section{Equivariant symbol map for bidifferential operators}
\label{sec.bidiffl}
Our method for building a projectively equivariant symbol map relies on 
the comparison of the spectra of $\casisymb$ and $\casi$.
For the sake of convenience, we shall first do this comparison, and 
show how it leads to a symbol map, under a technical hypothesis 
which we shall then weaken.

Let us name \emph{order} the index $i$ of an eigenvalue $\gamma_{i,p}$ of 
$\casisymb$.
Lemma~\ref{lem.eigen} shows that the shift $\delta$ is fixed by any 
equality of two eigenvalues with different orders.
\begin{defi}
   The value of the shift $\delta$ is said to be \emph{resonant} if 
   there exist $i,p,j,q\in\N$ satisfying $p\leq 
   \pent{i/2},q\leq\pent{j/2}$ and $i>j$ such that the two eigenvalues 
$\gamma_{i,p}$ and $\gamma_{j,q}$ of 
$\casisymb$ are equal.
\end{defi}
A resonant value will then be denoted by $\delta_{i,p;j,q}$ for a 
suitable choice of the indices.
\subsection{Symbol map when the shift is not resonant}
\begin{lem}
Assume that $\delta$ is not resonant.  Then
\begin{enumerate}
    \item  $\casisymb$ and $\casi$ have the same spectrum;

    \item  each eigenvector $P$ of $\casisymb$ belongs to some space 
    $\S{(i,p)}$;

    \item  such a $P$ is the principal symbol of a unique 
    eigenvector $Q(P)$ of $\casi$, which is associated to the same 
    eigenvalue.
\end{enumerate}
\end{lem}
\begin{proof}
   Let first $P \in \msymb{2}{\delta}$ be such that $\casisymb 
   P=\gamma_{i,p} P$.  
   Since $\delta$ is not resonant and $p\mapsto \gamma_{i,p}$ is 
   injective, $P\in \S{(i,p)}$.
   
   Assume now that $P$ is an eigenvector of $\casi$ of eigenvalue 
   $\gamma$ and that $P$ has total order $i$. 
   Let us rewrite the condition
   \begin{equation}\label{eq:eigen}                                                                                                                     
       \casi P = \gamma P                                                                                                                                   
   \end{equation}
   according to the decomposition $P=\sum_{j,q} P_{j,q}$ with 
   $P_{j,q}\in\S{(j,q)}$ and $P_{j}=\sum_{q} P_{j,q}$. We get
   \begin{equation}\label{eq:eigen-det}
       \left\{
       \begin{array}{llll}
	   \sum_{q\leq\pent{i/2}}(\gamma-\gamma_{i,q})P_{i,q} & = & 0 & \\
	   \sum_{q\leq\pent{j/2}}(\gamma - \gamma_{j,q})P_{j,q} & = & 
	   \casinil P_{j+1} & j=0,\ldots,i-1.
       \end{array}
       \right.
   \end{equation}
   The first equation shows that there exists one and only one $q_{0}$ 
   such that $\gamma=\gamma_{i,q_{0}}$, and that $P_{i}=P_{i,q_{0}}$ is an 
   eigenvector of $\casisymb$.  The second, that $P$ is uniquely determined by $P_{i,q_{0}}$. 
   
   We may thus define $Q(P_{i,q_{0}})=P$.  The value of the map $Q$ is clearly 
   well-defined on any eigenvector of $\casisymb$, hence the 
   proof.
\end{proof}

The correspondence $Q$ will now be considered as linearly extended to 
$\msymb{2}{\delta}$.
\begin{lem}
If $\delta$ is not resonant, then $Q$ is 
$\algproj$-equivariant.
\end{lem}
\begin{proof}
    Let $P\in\S{(i,p)}$. It suffices to prove that
    \[
    Q(\Lie{X}P)  = \Lie{X}(Q(P))
    \]
    for any $X\in\algproj$.  Since $\Lie{X}P$ also belongs to 
    $\S{(i,p)}$, both members of this equality are eigenvectors of 
    eigenvalue $\gamma_{i,p}$ of $\casi$.  They have the same principal 
    symbol and are thus equal. 
\end{proof}

\begin{thm}
If $\delta$ is not resonant, there exists a unique 
$\algproj$-equivariant 
symbol map on $\mdiff 2 \lambda \mu$.  It is a 
differential operator.
\end{thm}
\begin{proof}
    Taking account of the preceding lemmas, we just have to show that 
    the map $Q$ is such that $Q^{-1}$ is the unique equivariant symbol map and that it 
    is bijective.  
    
    On the one hand, any analogue of $Q$
    prolongs an eigenvector of $\casisymb$ into an 
    eigenvector of $\casi$ with the same eigenvalue and principal 
    symbols, and hence equals $Q$.
    
    Now, if $Q(P)=0$, the principal symbol of $P$, and therefore $P$, 
    equal $0$.  If $D\in\mdiff 2 \lambda \mu$ admits $P$ as its 
    principal symbol, then $Q(P)-D$ has an order less than that of 
    $D$ and an evident induction allows to conclude that $Q$ is 
    surjective.
\end{proof}

\subsection{Critical shift values}

A resonant value of the shift does not automatically prevent us from 
extending eigenvectors of $\casisymb$ into eigenvectors of $\casi$.  
Consider indeed $P_{i,p}\in\S{(i,p)}$ and assume that \mbox{$\gamma_{j, 
q}=\gamma_{i,p}$} for some $j<i$.  The equation (\ref{eq:eigen-det}) shows that 
one will still be able to prolong $P$ in an eigenvector of $\casi$ 
if the component of $\casinil P_{j+1}$ in $\S{(j,q)}$ vanishes.

We thus want to describe the spaces of homogeneous eigenvectors of 
$\casisymb$ actually reached through the application of $\casinil$ on 
one of them.

\begin{lem}\label{lem:reach}
    If $k,l,q\in\N$ are such that $q\leq\min(k,l)$,
\[
\casinil(\S{k,l,q})\subset\S{k-1,l,q-1}\oplus\S{k-1,l,q}\oplus\S{k,l-1,q-1}\oplus\S{k,l-1,q}.
\]
Hence, for $i,p\in\N$ such that $p\leq\pent{i/2}$,
\[
\casinil(\S{(i,p)})\subset \S{(i-1,p-1)}\oplus\S{(i-1,p)}.
\]
\end{lem}
\begin{proof}                                                                                                                                        
    Notice first 
    that                                                                                                                                    
    \begin{eqnarray*}                                                                                                                                    
	(\eta D_\alpha) \hwv &=& q \eta_1 (\alpha_1\beta_2 - 
	\alpha_2\beta_1)^{q-1} \alpha_1^{k-q}\beta_1^{l-q}\beta_2\\
	&&\mbox{}- q \eta_2 (\alpha_1\beta_2 - 
	\alpha_2\beta_1)^{q-1}\alpha_1^{k-q}\beta_1^{l-q+1}\\                                                         
	&&\mbox{}+(k-q) \eta_1 (\alpha_1\beta_2 - 
	\alpha_2\beta_1)^q\alpha_1^{k-q-1}\beta_1^{l-q}.                                                           
    \end{eqnarray*}                                                                                                                                      
    The second and third terms are respectively multiples of 
    $\hwv[k-1,l,q-1]$ and
    $\hwv[k-1,l,q]$, and the first 
    equals
    \[                                                                                                                                                   
    \frac{1}{k+l+1-2q}(\rho(e_2\otimes 
    \epsilon^1)\hwv[k-1,l,q-1]-(k-q)\hwv[k-1,l,q]),                                                                  
    \]                                                                                                                                                   
    where 
    $\rho(e_2\otimes\epsilon^1) = \alpha_2 D_{\alpha_1}+\beta_2 
    D_{\beta_1}$.             
    Therefore,                                                                                                                                           
    \[                                                                                                                                                   
    (\eta D_\alpha)\hwv \in 
    \S{k-1,l,q-1}\oplus\S{k-1,l,q}.                                                                                              
    \]                                                                                                                                                   
    Let 
    $P\in\S{k,l,q}$.                                                                                                                                        
    As                                                                                                                                                   
    \[                                                                                                                                                   
    (\eta D_\alpha) \rho(h\otimes\zeta) P = <h,\eta>(\zeta D_\alpha) P 
    + \rho(h\otimes\zeta)(\eta D_\alpha) 
    P,                  
    \]                                                                                                                                                   
    we observe that $(\eta D_\alpha) P$ also belongs to 
    $\S{k-1,l,q-1}\oplus\S{k-1,l,q}$.                                                                
    The result then follows from the definitions of $\casinil$ and 
    $\S{(i,p)}$.
\end{proof}   

Remember that the value of the shift $\delta$ is \emph{resonant} if 
   there exist $i,p,j,q\in\N$ satisfying the following conditions:
   \begin{enumerate}
      \item $p\leq \pent{i/2}$\label{eq.condcritp},
      \item  $q\leq\pent{j/2}$\label{eq.condritq},
      \item $i>j$\label{eq.condcriti},
      \item  $\gamma_{i,p}=\gamma_{j,q}$\label{eq.condcritg}.
   \end{enumerate}

Lemma~\ref{lem:reach} leads us to the following definition.
\begin{defi}
    A resonant value $\delta$ is \emph{critical} if it  may be writte 
    $\delta_{i,p;j,q}$ with $i,p,j,q$ satisfying
\begin{enumerate}
 \setcounter{enumi}{4}
    \item $0\leq p-q\leq i-j$\label{eq.condcrit}
\end{enumerate}
in addition to conditions~(\ref{eq.condcritp})--(\ref{eq.condcritg}) 
above.
\end{defi}

Before proving existence and uniqueness of a projectively equivariant 
symbolization when $\delta$ is resonant but not critical, we shall 
examine  these values somewhat more precisely.

A shift of $0$ occurs  in the most ``natural'' uses of densities~: 
when the weight of each considered density is $0$, 
which                           
means the multidifferential operators act on functions, and when 
$\mu=\sum_i \lambda_i$, which is required to consider 
the                           
natural product of densities.  
But $0$ is a resonant value of $\delta$, whatever the dimension $n>1$. 
Notice indeed that if  $\delta=0$ and $k=n-1$, then $\gamma_{6+k,3}=\gamma_{5+k,0}$.

These remarks explain why we wanted to prove that 
zero 
is not a critical value of the shift. One has the following stronger 
property.                                                                                                                           
                                                                                                                                                   
\begin{prop}\label{prop.critshift}                                                                                                                                  
   Let $f:\N\setminus\{0\}\rightarrow \Q:i\mapsto 
   \delta_{i,\pent{i/2};0,0}$.
    Then
    \begin{itemize}
	\item $f(1)=1$,
	\item $f$ is increasing and
	\item $\lim_{i\rightarrow+\infty}f(i)=+\infty$.
    \end{itemize}
    Moreover, if $i,j,p,q\in\N$ satisfy    
    conditions~(\ref{eq.condcritp})--(\ref{eq.condcrit}) above,
    then 
    \[\delta_{i,p;j,q}\geq f(i).\]
\end{prop}                                                                                                                                    
\begin{proof}                                                                                                                                        
   The first three assertions are evident since 
   \begin{eqnarray*}
       \delta_{2k,k;0,0}&=&1+\frac{3(k-1)}{2(n+1)}, \quad\forall k\in\N\setminus\{ 
       0\}\\  
   \text{and }\delta_{2k+1,k;0,0}&=&1+\frac{3k^2}{(2k+1)(n+1)}, \quad\forall k\in\N.
   \end{eqnarray*}
    
   To prove the stated inequality, observe first that $\delta_{i,p,j,q}>\delta_{i,p',j,q}$ if 
    $p<p'\leq \pent{i/2}$.  Indeed, one has then
    \[
    \delta_{i,p;j;q}-\delta_{i,p+1;j,q}=\frac{i-2p}{i-j}>0.
    \]
     
    Now, if $\pent{i/2}>i-j+q$ then
	\begin{eqnarray*}
	\delta_{i,p;j,q}&\geq& \delta_{i,i-j+q,j,q}\\
	&=&\delta_{i,i-j+q;0,0}+\frac{1}{(n+1)i}\overbrace{(j(i-j-1)+q(2j-q)+q)}^{\geq 
	0}\\
	&\geq& \delta_{i,i-j+q;0,0}\\
	&\geq& \delta_{i,\pent{i/2};0,0}.
	\end{eqnarray*}
    
	Else, suppose that $i=2k+1$ and notice that $k+1-j+q\geq 0$.
    One can write
    \[
    \delta_{i,p;j,q}\geq \delta_{2k+1,k;j,q}.
    \]
    Moreover,
    \begin{multline*}
    (\delta_{2k+1,k;j,q}-\delta_{2k+1,k;0,0})(2k+1)(2k+1-j)(n+1)=\\
    3k^2 
    j+2kjq-2kj^2-2kq^2
    +2kj+2jq-j^2+2kq-q^2
    +j+q.
    \end{multline*}
    Denote by $a_{r}(k,j,q)$ the sum of the homogeneous monomials of 
    degree $r$ in the latter expression. It is clear that 
    $a_{2}(k,j,q)\geq 0$.  
    Furthermore, if $q\leq k-2$,
    \begin{eqnarray*}
	a_{3}(k,j,q)+j&=& k((k-2)j+2(k+1-j+q)j-2q^2)+j \\
	&\geq&k((k-2)j-2q^2)+j\\
	&\geq&j(k(k-2-q)+1)\\
	&\geq&0.
    \end{eqnarray*}
    Since $2q\leq j\leq 2k$, it suffices to compute explicitly the different values 
    of $a_{3}(k,j,q)+j$ for $k-1\leq q\leq k$ to prove that 
    $\delta_{i,p;j,q}\geq \delta_{i,\pent{i/2};0,0}$ for any odd 
    natural number $i$.
   
    One proceeds in the same way if $i=2k$ and $k-j+q\geq 0$.
\end{proof}                                                                                                                                          

\begin{cor}
    The number of critical values of the shift $\delta$ contained in 
    any interval $[a,b]$ is finite and vanishes if $b<1$.
\end{cor}
We shall eventually show~(see Section~\ref{subsec.crit}) that $1$  is a 
critical value, for which there exists an equivariant symbol map if 
and only if $\lambda=(0,0)$.  Furthermore, it suffices to 
examine the case of second order operators to prove this.

\subsection{Symbol map for resonant values of $\delta$}
Let now $\delta$ be resonant but not critical and let $P\in\S{(i,p)}$.  We can 
still define $Q(P)$, the unique eigenvector of $\casi$ admitting $P$ as 
its principal symbol and belonging to 
\[
\widetilde{\S{(i,p)}}=\bigoplus_{j,q\,:\,0\leq p-q\leq i-j} \S{(j,q)}.
\]
\begin{prop}\label{pro:equiv2}
    The map $Q$,  linearly extended to $\msymb{2}{\delta}$,  
 is equivariant.
\end{prop}

\begin{lem}\label{lem.equiv}
    Let $i,p\in\N$ be such that $p\leq\pent{i/2}$ and $X\in\algproj$.
    If the polynomial form of an operator $T$ belongs to 
    $\widetilde{\S{(i,p)}}$, then so does the polynomial form of 
    $\Lie{X} T$.
\end{lem}
\begin{proof}[Proof of lemma \ref{lem.equiv}]
    Assume that $P\in\S{k,l,p}$, $(k+l=i)$, is the polynomial form of $T$.  
    As it can be seen from expressions~(\ref{eq.LieDerT}) 
    and~(\ref{eq.LieDerP}), the polynomial form of $\Lie{X}T$ differs 
    from $\Lie{X}P$ by
    \begin{multline*}
	 - \lambda_{1}<X,\zeta> (\zeta D_{\alpha})P 
	-\lambda_{2}<X,\zeta> (\zeta D_{\beta})P \\
	-\frac{1}{2}<X,\alpha> (\zeta D_{\alpha})^2 P- 
	\frac{1}{2}<X,\beta>(\zeta D_{\beta})^2 P.
    \end{multline*}    
   This quantity vanishes if the degree of $X$ does not exceed 
   one.  We have shown in the proof of Lemma \ref{lem:reach} that both 
   $<X,\zeta>(\zeta D_{\alpha})P$ and $<X,\zeta>(\zeta D_{\beta})P$ 
   belong to $\widetilde{\S{(i,p)}}$ for any vector field 
   $X\in\algproj$.
     When $X=\theta^*$, the last two terms equal $-k (\theta D_{\alpha})P-l(\theta D_{\beta})P$ and 
    thus also belong to $\widetilde{\S{(i,p)}}$.
\end{proof}
\begin{proof}[Proof of proposition \ref{pro:equiv2}]
    Proceed as when $\delta$ is not resonant, noticing that both 
    members of the equality
    \[
    \Lie{X}Q(P)=Q(\Lie{X}P)
    \]
    belong to $\widetilde{\S{(i,p)}}$ provided that $P\in\S{(i,p)}$.
\end{proof}

To prove the uniqueness of the symbol map, we must ensure that 
the prolongation of $P\in\S{(i,p)}$ into an eigenvector of $\casi$ does not get out of 
$\widetilde{\S{(i,p)}}$.  In the statement of the next proposition, 
notice that tensor fields and bidifferential  
operators are isomorphic modules over the constant and linear vector fields.

\begin{prop}
   Any linear map $T : \msymb{2}{\delta}\rightarrow \mdiff 2 \lambda \mu$ that is  
    equivariant with respect to the constant and linear vector fields stabilizes each space $\widetilde{\S{(i,p)}}$. 
\end{prop}
\begin{proof}
    A slight adaptation of the proof of theorem 5.1 of \cite{lo} 
    shows that  such a map is local. Being invariant under the action 
    of constant vector fields,  
    it is thus a differential operator with constant coefficients.  The invariance 
    of $T$ with respect to the linear vector fields means that the polynomial 
    form
    of $T$ maps
    \[
    X^r\otimes Y^s \quad (X,Y\in\R^n)
    \]
    onto a tensor field with homogeneous component in 
    $\vee^p T\R^n\otimes\vee^q T\R^n$ given by
    \begin{eqnarray*}
        (\alpha,\beta)&\mapsto &\sum_{r_{1},r_{2}} c_{r_{1},r_{2}}  
    <X,\alpha>^{r_{1}}<X,\beta>^{r_{2}}<X,\eta>^{r-(r_{1}+r_{2})}\\
    &&\mbox{}\times<Y,\alpha>^{p-r_{1}}<Y,\beta>^{q-r_{2}}<Y,\eta>^{s-(p+q)+(r_{1}+r_{2})}
    \end{eqnarray*}
     where $\eta$ symbolizes the derivatives of the argument of $T$.
    It means that the polynomial symbolization of $T$ is a linear 
    combination of powers of the operators $\alpha D_{\beta}, \beta 
    D_{\alpha}, \eta D_{\alpha}$ and $\eta D_{\beta}$.
    
    But the last two stabilize the space 
    $\widetilde{\S{(i,p)}}$, because of the definition of this space.  
    Furthermore, $\alpha D_{\beta}$ and $\beta D_{\alpha}$ vanish or
    intertwine irreducible $\alg{sl}{n,\R}$-submodules of the spaces
    $\S{(i,p)}$.  Hence the conclusion.
\end{proof}
\subsection{The one-dimensional case}
\label{sec.n1}
Assume that $n=1$.  The formula~(\ref{eq.CasiDef}) still holds.  One sees easily 
that the spaces
\[
\bigoplus_{k+l=i} \Gamma(\vee^k T\R\otimes \vee^l 
T\R\otimes\Delta^\delta \R)
\]
are made of homogeneous eigenvectors of $\casisymb$ which are 
associated to the eigenvalue $\gamma_{i,0}$ given 
by~Lemma~\ref{lem.eigen} when $n=1$ and $q=0$.

The symbol map may be built as before when the shift is not resonant 
(i.e. when no eigenvalues $\gamma_{i,0}$ and $\gamma_{j,0}$, $(i\neq 
j)$, are equal). 
Moreover, $\gamma_{i,0}=\gamma_{j,0}$ if and only if 
$\delta=1+\frac{i+j-1}{2}$.  These are the resonant values of the 
shift.  They are all critical in the sense that iterated applications 
of $\casinil$ on $\S{(i,0)}$ always end up in $\S{(j,0)}$ (cf.~Lemma~\ref{lem:reach}).

\section{Second order operators}
\label{sec.secord}
Let us now give an explicit equivariant symbol map example by carrying some 
computations for second order operators.  We shall also observe that 
a critical shift value may prevent the symbol map from existing or 
being 
unique.

Examining the natural action of $\alg{sl}{n,\R}$ on highest weight 
vectors described in formula~(\ref{eq:hwv}), one easily checks that
\[
    \S{(2,0)}=\S{2,0,0}\oplus\S{0,2,0}\oplus\S{1,1,0}
\]
and that
\[
\S{(2,1)}=\S{1,1,1}.
\]
This means that elements of $\S{(2,0)}$ are symmetric quadratic polynomials 
of $\alpha$ and $\beta$ while $\S{(2,1)}$ contains the antisymmetric 
ones.

Shift values for which an equality between the eigenvalues 
$\gamma_{2,0}$, $\gamma_{2,1}$, $\gamma_{1,0}$ and 
$\gamma_{0,0}$ of $\casisymb$ may occur are given in 
Figure~\ref{fig.crit}. Notice that they are all critical.

\begin{figure}
\begin{center}
    \vspace{1em}
\begin{tabular}{|cc|c|c|}\hline
 &$\gamma_{2,0}=\gamma_{1,0}$&$\gamma_{2,0}=\gamma_{0,0}$&$\gamma_{2,1}=\gamma_{1,0}=\gamma_{0,0}$
 \\[.5ex] \hline
&$\delta=\frac{n+3}{n+1}$&$\delta=\frac{n+2}{n+1}$&$\delta=1$\\[.5ex]\hline
\end{tabular}
\vspace{1em}
\end{center}
\caption{Resonant (and critical) shift values for second order operators}
\label{fig.crit}
\end{figure}

\subsection{Generic shift values}
Assume that $\delta$ takes none of the above listed values.
Using notations of the equations~(\ref{eq:eigen-det}), we  
want first to associate 
to the monomial $P_{2}=\alpha_{i}\alpha_{j}$  lower degree tensors 
$P_{1}$ and $P_{0}$ such that $P_{2}+P_{1}+P_{0}$, considered as a 
bidifferential operator, is an eigenvector of $\casi$.
 The sought terms
$P_{1}$ and $P_{0}$ are the only solutions of the system
\begin{equation}\label{eq:eigen-det2}
\left\{
\begin{array}{lcl}
    (\gamma_{2,0}-\gamma_{1,0})P_{1} &=& 2((n+1)\lambda_{1}+1)(\eta 
    D_{\alpha})P_{2} \\
    (\gamma_{2,0}-\gamma_{0,0})P_{0} &=& 2(n+1)\lambda_{1}(\eta 
    D_{\alpha})P_{1}
\end{array}
\right.
\end{equation}
 and the inverse of the symbol map will thus associate to the tensor field 
\[
c\, \alpha_{i}\alpha_{j},
\]
where $c$ belongs to $\Gamma(\Delta^\delta\R^n)$, 
the operator
\begin{multline}\label{eq:eigen-det2a}
    (f,g)\mapsto c\,\del_{i}\del_{j}f \, g \\
    \mbox{}+ 
\frac{(n+1)\lambda_{1}+1}{(n+1)(1-\delta)+2}\,(\del_{i}c\,\del_{j}f\, g+\del_{j}c\,\del_{i}f\, g)\\
\mbox{}+\frac{((n+1)\lambda_{1}+1)(n+1)\lambda_{1}}{((n+1)(1-\delta)+2)((n+1)(1-\delta)+1)}\,\del_{i}\del_{j}c\, f\, g.
\end{multline}

Exchanging $f$ and $g$ and substituting $\lambda_{2}$ to $\lambda_{1}$ 
in the last operator gives the one to be associated to $c\,
\beta_{i}\beta_{j}$.  
Using the same notations, it is easy to compute that one should associate to the tensor field
\[
c\, (\alpha_{i}\beta_{j}+\alpha_{j}\beta_{i})
\]
the operator
\begin{multline}\label{eq:eigen-det2ab}
    (f,g)\mapsto  c\, (\del_{i} f \del_{j} g+ \del_{j} f \del_{i} g)\\
    \mbox{}+ \frac{(n+1)\lambda_{1}}{(n+1)(1-\delta )+2}
    (\del_{i} c \, f\del_{j}g + \del_{j} c \, f \del_{i} g)\\
  \mbox{}+  \frac{(n+1)\lambda_{2}}{(n+1)(1-\delta )+2}
    (\del_{j} c \, \del_{i}f g + \del_{i} c \,  \del_{j} f g)\\
    \mbox{}+2  \frac{(n+1)^2\lambda_{1}\lambda_{2}}
    {((n+1)(1-\delta)+2)((n+1)(1-\delta)+1)}\del_{i}\del_{j}c\, f g,
\end{multline}
and 
\begin{multline}\label{eq:eigen-det2As}
   (f,g) \mapsto  c\, (\del_{i}f\del_{j}g-\del_{j}f\del_{i}g)\\
  \mbox{}+\frac{\lambda_{1}}{1-\delta}(\del_{i}c\, 
  f\del_{j}g-\del_{j}c\,f\del_{i}g)
+\frac{\lambda_{2}}{1-\delta}(\del_{j}c\,\del_{i}f g-\del_{i}c\,\del_{j}f g)
\end{multline}
to the tensor field
\[
c \,(\alpha_{i}\beta_{j}-\alpha_{j}\beta_{i}).
\]
To degree 1 tensors
\[
c \,\alpha_{i},
\]
one associates

\begin{equation}\label{eq:eigen-det1}
(f,g)\mapsto c\, \del_{i}f g+ \frac{\lambda_{1}}{1-\delta}\del_{i}c f g
\end{equation}
and similarly, substituting $\lambda_{2}$ to $\lambda_{1}$, to degree 1 tensors in $\beta$.

Notice that tensors of degree 
$0$ and  bidifferential operators of order $0$ make up isomorphic 
$\vect(\R^n)$-modules.

For each generic value of $\delta$, we have thus described a linear bijection between
 $\S{(2,0)}\oplus\S{(2,1)}\oplus\S{(1,0)}\oplus\S{(0,0)}$ and the 
submodule of $\mdiff 2 \lambda \mu$ made up of operators with total 
order not greater than $2$. This bijection is the only one to be 
projectively equivariant and to preserve the principal symbol of its 
arguments.
\subsection{Critical shift values}
\label{subsec.crit}
Assume that $\delta = \frac{n+3}{n+1}$. Since 
$\gamma_{2,0}=\gamma_{1,0}$, the first  equation 
of~(\ref{eq:eigen-det2}) may be satisfied only if 
$\lambda_1=\frac{-1}{n+1}$ ; when rewritten for eigenvectors of the 
form $\beta_{i}\beta_{j}$, this equation forces 
$\lambda_{2}=\frac{-1}{n+1}$.

Now, it may be seen from expression~(\ref{eq:eigen-det2ab}) that no 
eigenvector of $\casi$ admits $\alpha_{i}\beta_{j}+\alpha_{j}\beta_{i}$ 
as its principal symbol, since this would require $\lambda_{1}$ or 
$\lambda_{2}$ to vanish.

We have just proved the following statement.
\begin{prop}
    For all values of $\lambda$ and $\mu$ such that 
    $\mu-\lambda_{1}-\lambda_{2}=\frac{n+3}{n+1}$, there  exists no 
    $\algproj$-equivariant symbol defined on $\mdiff 2 \lambda \mu$.
\end{prop}

\begin{prop}
    If $\delta=\frac{n+2}{n+1}$,
there exists a projectively equivariant symbol  for bidifferential 
second order operators of 
$\mdiff 2 \lambda \mu$ if and only if
\[
\lambda\in\{(0,\frac{-1}{n+1}),(\frac{-1}{n+1},0),(0,0)\},
\] in which 
case this symbol is not unique.
\end{prop}
\begin{proof}
 As it was done in the case $\delta=\frac{n+3}{n+1}$, it can be checked that one 
of the weights $\lambda_{1}$ and $\lambda_{2}$ must vanish and the 
other be chosen in $\{0,\frac{-1}{n+1}\}$.
A direct computation then allows to prove that the map 
which associates to $c\, \alpha_{i}\alpha_{j}\in \S{(2,0)}$ the 
bidifferential operator 
\begin{eqnarray*}
    (f,g)&\mapsto &c\,\del_{i}\del_{j}f \, g \\
    &&\mbox{}+ 
\frac{(n+1)\lambda_{1}+1}{(n+1)(1-\delta)+2}\,(\del_{i}c\,\del_{j}f\, g+\del_{j}c\,\del_{i}f\, g)\\
&&\mbox{}+k\,\del_{i}\del_{j}c\, f\, g
\end{eqnarray*}
is projectively equivariant for any value of the parameter $k$. One 
defines similarly an operator associated to $c\,\beta_{i}\beta_{j}$.
This equivariance also holds for the map that associates to 
$c(\alpha_{i}\beta_{j}+\alpha_{j}\beta_{i})\in\S{(2,0)}$ the operator
\begin{multline*}
    (f,g)\mapsto  c\, (\del_{i} f \del_{j} g+ \del_{j} f \del_{i} 
    g)\\
    \mbox{}+ (n+1)\lambda_{1}
    (\del_{i} c \, f\del_{j}g + \del_{j} c \, f \del_{i} g)
    +  (n+1)\lambda_{2}
    (\del_{j} c \, \del_{i}f g + \del_{i} c \,  \del_{j} f g)\\
    \mbox{}+k\,   
    \del_{i}\del_{j}c\, f g,
\end{multline*}
whatever the value of $k$.
Moreover, any eigenvector of $\casisymb$ not belonging to $\S{(2,0)}$ 
is the principal symbol of a unique eigenvector of $\casi$.
\end{proof}

Finally, when $\delta =1$, expression~(\ref{eq:eigen-det2As}) 
 shows that both $\lambda_{1}$ and $\lambda_{2}$ must vanish.
 Then, both maps 
\begin{equation}\label{eq:vectinv1}
T_{1}:c(\alpha_{i}\beta_{j}-\alpha_{j}\beta_{i})\in\S{(2,1)}\mapsto 
((f,g)\mapsto \del_{i}c f\del_{j}g-\del_{j} c f\del_{i}g)
\end{equation}
and
\begin{equation}\label{eq:vectinv2}
T_{2}:c\,\alpha_{i}\mapsto ((f,g)\mapsto \del_{i}c f g)
\end{equation}
are $\algproj$-equivariant~(cf. below).  Similar maps may be defined, switching 
the first and second arguments.
But we may release the requirement on the 
operators not to have an order greater than two and state the 
following proposition.
\begin{prop}
   If $\delta=1$, there exists a projectively equivariant
   symbol  for the elements of $\mdiff 2 \lambda \mu$ if and only if
$\lambda=(0,0)$, in which case this symbol is not unique.
\end{prop}
\begin{proof}
    Equation~(\ref{eq:eigen-det}) shows that the prolongation of any eigenvector 
    $P\in\S{(i,p)}$, $(i>2)$, of $\casisymb$ will be impossible only if 
    $\gamma_{i,p}=\gamma_{j,q}$ for some $j,q$ such that $i,p,j$ and 
    $q$ define a critical value of the shift. But then, in view of 
    Proposition~\ref{prop.critshift},
    $\delta=\delta_{i,p;j,q}\geq 
    \delta_{i,\pent{i/2};0,0}\geq\delta_{3,1;0,0}=\frac{n+2}{n+1}$, 
    hence a contradiction.
\end{proof}
\subsection{A link with the linear case}
We first remark that one can obtain formulas~(\ref{eq:eigen-det2a}),%
~(\ref{eq:eigen-det2ab}) and~(\ref{eq:eigen-det1}) from the following result,
which is an adaptation of~\cite[Prop. 4.4]{lo}.
\begin{prop}
    Let $\lambda,\mu\in\R$ be such that 
    $\mu-\lambda\not\in\{1,\frac{n+2}{n+1},\frac{n+3}{n+1}\}$.  Denote 
    by $\mathcal{D}^2_{\lambda,\mu}$ the $\vect(\R^n)$-module of 
    linear second order differential operators acting on 
    $\lambda$-densities and valued in $\mu$-densities. Let also
    $\mathcal{S}^{(1)}_{j}$ denote the module 
    $\Gamma(\vee^j 
    T\R^n\otimes\Delta^{\mu-\lambda}\R^n)$.
    Then there exists a unique isomorphism of $\algproj$-modules
    \[
    q_{\lambda,\mu}:
    \bigoplus_{0\leq j\leq 2}\mathcal{S}^{(1)}_{j}
    \rightarrow \mathcal{D}^2_{\lambda,\mu}
    \]
    that preserves the principal symbol of its arguments.
\end{prop}
Let us define
\begin{eqnarray*}
\tau_{\alpha}&:&\mathcal{S}^{(1)}_{2}\rightarrow\S{2,0,0}:P\mapsto((\alpha,\beta)\mapsto P(\alpha)),\\
\tau_{\beta}&:&\mathcal{S}^{(1)}_{2}\rightarrow\S{0,2,0}:P\mapsto((\alpha,\beta)\mapsto P(\beta))\\
\text{and }\tau_{\alpha\beta}&:&\mathcal{S}^{(1)}_{2}\rightarrow\S{1,1,0}:P\mapsto((\alpha,\beta)\mapsto P(\alpha+\beta)-P(\alpha)-P(\beta)).
\end{eqnarray*}
These are  isomorphisms of $\vect(\R^n)$-modules.
The operator~(\ref{eq:eigen-det2a}) is then nothing but 
\[
(q_{\lambda_{1},\mu-\lambda_{2}}\circ\tau_{\alpha}^{-1}(c \,
\alpha_{i}\alpha_{j}))f\cdot g.
\]
One  obtains in a like manner expression~(\ref{eq:eigen-det1}).
Similarly, operator~(\ref{eq:eigen-det2ab}) is associated to the 
polynomial $P=c\,(\alpha_{i}\beta_{j}+\alpha_{j}\beta_{i})$ by the 
formula
\begin{multline*}
(q_{\lambda_{1}+\lambda_{2},\mu}\circ\tau_{\alpha\beta}^{-1}(P))(fg)\\
\mbox{}-(q_{\lambda_{1},\mu-\lambda_{2}}\circ\tau_{\alpha}^{-1}(c \,
\alpha_{i}\alpha_{j}))f\cdot g-f\cdot 
(q_{\lambda_{2},\mu-\lambda_{1}}\circ\tau_{\beta}^{-1}(c \,
\beta_{i}\beta_{j}))g.
\end{multline*}

\subsection{Equivariant maps}
Our second remark regards the maps $T_1$ and $T_2$ defined in~(\ref{eq:vectinv1}) 
and~(\ref{eq:vectinv2}).  They are in fact $\vect(\R^n)$-equivariant 
and can be interpreted as follows.

On any oriented manifold $M$ equipped with a nowhere vanishing volume 
form $\omega$, densities of weight $1$ and volume forms make up 
isomorphic modules over $\vect(M)$.  Any element of $\Gamma(\bigwedge^2 
\mathit{TM}\otimes \Delta^1 M)$ --- the analogue of our space $\S{(2,1)}$ --- 
may be written $\Lambda \,\omega$ 
with $\Lambda\in\bigwedge^2 \mathit{TM}$.
Then in any chart domain where $\omega$ admits the local form 
$dx^1\wedge\cdots\wedge dx^n$, $T_1$ is the local form of the globally 
defined operator
\[
\Lambda\omega\mapsto ((f,g)\mapsto f\cdot dg\wedge 
d(i(\Lambda)\omega))
\]
where $f$ and $g$ are smooth functions.
Similarly, $T_2$ is the local form of 
\[
X\omega\mapsto((f,g)\mapsto L_{X}\omega \cdot f\cdot g)
\]
with $X\in\vect(M)$.

\section*{Acknowledgements}
We would like to thank M. De Wilde, P. Lecomte and P. Mathonet for 
helpful suggestions and V. Ovsienko for his interest in this work, 
particularly during our stay at the C.P.T.-C.N.R.S. in Luminy.

This work was supported by a Research Fellowship of the Belgian National Fund for 
    Scientific Research.

\end{document}